\documentclass[11pt]{article}
\usepackage{amsmath,amssymb,latexsym}
\usepackage{theorem}
\newtheorem{theorem}{Theorem}
\newtheorem{lemma}{Lemma}

\newtheorem{proposition}{Proposition}
{\theorembodyfont{\rmfamily} }

\newcommand{\sen}{\mathop{\rm sen}}

\newcommand{\supp}{\mathop{\rm supp}}
\newcommand{\field}[1]{\mathbb{#1}}
\newcommand{\R}{\field{R}}

\newcommand{\C}{\field{C}}
\newcommand{\T}{\field{T}}
\newcommand{\D}{\field{D}}

\title{Sobolev orthogonal polynomials: balance and asymptotics }
\author{Manuel Alfaro\thanks{Partially supported by MEC of Spain under
Grant MTM2006-13000-C03-03, FEDER funds (EU), and the DGA project
E-64 (Spain).}
\\ Departamento de Matem\'{a}ticas and IUMA. Univ. de Zaragoza (Spain)
\and Juan Jos\'{e} Moreno--Balc\'{a}zar \thanks{Partially supported by MEC
of Spain under Grant MTM2005--08648--C02--01 and Junta de Andaluc\'{\i}a
(FQM229 and excellence projects FQM481, P06-FQM-1735).}
 \\ Departamento de Estad\'{\i}stica y Matem\'{a}tica Aplicada \\ Univ. de Almer\'{\i}a (Spain)
\\ Instituto Carlos I de F\'{\i}sica Te\'orica y Computacional \\
Univ. de Granada (Spain) \and  Ana Pe\~{n}a \thanks{Partially supported
by MEC of Spain under Grants MTM 2004-03036 and
MTM2006-13000-C03-03,
FEDER funds (EU), and the DGA project E-64 (Spain).} \\
Departamento de Matem\'{a}ticas. Univ. de Zaragoza (Spain) \and M. Luisa
Rezola$^{*}$ \\ Departamento de Matem\'{a}ticas and IUMA. Univ. de
Zaragoza (Spain) }

\date{}

\begin{document}

\maketitle

\begin{abstract}
Let $\mu_0$ and $\mu_1$  be measures supported on an unbounded
interval and $S_{n,\lambda_n}$ the extremal varying Sobolev
polynomial which minimizes
\begin{equation*}
\langle P, P \rangle_{\lambda_n}=\int P^2 \, d\mu_0 + \lambda_n \int P'^2 \,
d\mu_1, \quad \lambda_n >0
\end{equation*}
\noindent in the class of all monic polynomials of degree $n$. The
goal of this paper is twofold. On one hand, we discuss how to
balance both terms of this inner product, that is, how to choose a
sequence $(\lambda_n)$ such that both measures $\mu_0$ and $\mu_1$
play a role in the asymptotics of $\left(S_{n, \lambda_n} \right).$
On the other, we apply such ideas to the case when both $\mu_0$ and
$\mu_1$ are Freud weights. Asymptotics for the corresponding $S_{n,
\lambda_n}$ are computed, illustrating the accuracy of the choice of
$\lambda_n\, .$
\end{abstract}



\bigskip
Corresponding author:

M. L. Rezola

Departamento de Matem\'aticas

Universidad de Zaragoza

50009-ZARAGOZA  (SPAIN)

e-mail: rezola@unizar.es

FAX: (+34) 976761338

\newpage

\section{Introduction}
\setcounter{equation}{0}

One of the central problems in the analytic theory of orthogonal
polynomials is the study of their asymptotic behavior. In this paper
we are concerned with the asymptotic properties of Sobolev
orthogonal polynomials, that is polynomials orthogonal with respect
to an inner product involving derivatives. In this sense, given
$\mu_0$ and $\mu_1$ finite Borel measures supported on an interval
$I \subset \R$ and $\lambda >0$ we consider the Sobolev inner
product
\begin{equation}\label{producto-lambda}
\langle  P,Q \rangle _{\lambda}=\int P\,Q\,d\mu_0+\lambda\,\int P'\,Q'\,d\mu_1
\end{equation}
in the space of all polynomials with real coefficients.

We denote by $P_{n,\mu_0}$, $P_{n,\mu_1}$ and $S_{n,\lambda}$ the
corresponding monic polynomials orthogonal with respect to $\mu_0$,
$\mu_1$ and $\langle  \cdot, \cdot  \rangle_{\lambda}$,
respectively.

Let $\mu_0$ and  $\mu_1$ be measures compactly supported on $\R$.
Whether $(\mu_0, \mu_1)$ is a coherent pair, which means that there
exist nonzero constants $\sigma_n$ such that the corresponding monic
polynomials satisfy for each $n$
\begin{equation*}
P_{n,\mu_1} = \displaystyle \frac{P'_{n+1,\mu_0}}{n+1} + \sigma_n
\frac{P'_{n,\mu_0}}{n}
\end{equation*}
or, if $\mu_0$ and  $\mu_1$ fulfill much milder conditions, i.e.,
they belong to the well known Szeg\H{o} class, it has been
established (see \cite{MMPP} and \cite{martinez:b}) that the ratio
asymptotics
\begin{equation*}
\lim_{n \to \infty} \frac{S_{n,
\lambda}(z)}{P_{n,\mu_1}(z)}=\frac{2}{\varphi'(z)}
\end{equation*}
holds uniformly on compact subsets of $\overline{\C} \setminus
[-1,1]$, where $\varphi (z)=z+\sqrt{z^2-1}$ with $\sqrt{z^2-1}>0$
when $z>1$. In other words, the measure $\mu_0$ does not appear
explicitly within the asymptotic expression.

Nevertheless, a closer look at the inner product
(\ref{producto-lambda}) explains the \lq \lq dominance" of the
measure $\mu_1$ in the asymptotics: the derivative makes the leading
coefficient of the polynomials in the second integral of
(\ref{producto-lambda})  be multiplied by the degree of the
polynomial. Thus, if we want both measures to have an impact on the
behavior of the polynomials for $n\to\infty$, it seems natural to
\lq\lq balance" the inner product, that is, to compensate both
integrals by introducing a varying parameter $\lambda_n$.

In a general framework, we consider the varying Sobolev inner
product $\langle  P,Q  \rangle_{{\lambda}_n}.$ We denote by
$S_{n,{\lambda}_n}$ the monic polynomial which minimizes the
expression $\langle  Q_n,Q_n  \rangle_{{\lambda}_n}$ in the class
of all monic polynomials $Q_n$ of degree $n$.

Concerning the choice of the varying parameter $\lambda_n$, it is
interesting to write the expression of the Sobolev inner product in
terms of monic polynomials, that is
\begin{equation}\label{normaQ_n-acotado}
\langle Q_n,Q_n \rangle_{\lambda_n}=\int
(Q_n)^2\,d\mu_0+\lambda_n\,n^2\int\left(\frac{Q'_n}{n}
\right)^2\,d\mu_1.
\end{equation}
In this expression each integral in the right hand side is bounded
from below by $\int P_{n,\mu_0}^2\,d\mu_0$ and $\int P_{n-1,
\mu_1}^2\,d\mu_1,$ respectively, as long as $Q_n$ is a monic
polynomial of degree $n$.

If the measures $\mu_0$ and $\mu_1$ are supported on the same
bounded interval where they satisfy the Szeg\H{o} condition, then
$\int P_{n,\mu_0}^2\,d\mu_0$ behaves as $\int P_{n-1,\mu_1}^2\,
d\mu_1$, when $n\to\infty$. More precisely, the ratio
$\displaystyle\frac{\int P_{n, \mu_0}^2\,d\mu_0}{\int P_{n-1,
\mu_1}^2\,d\mu_1}\,$ has a limit. Therefore, in order to balance
both terms in (\ref{normaQ_n-acotado}) it is natural to keep
$\lambda_n\,n^2$ bounded.

In fact, it was proved in \cite{AMR} that if $({\lambda}_n)$ is a
decreasing sequence of positive real numbers such that $\lim_{n}
{\lambda}_n\,n^2 \in (0, +\infty)$ then
$$\lim_{n \to \infty} \frac{S_{n,{\lambda}_n}(z)}{R_n(z)}=1$$
locally uniformly in $\overline{\C} \setminus [-1,1]$, where $(R_n)$
is the sequence of monic polynomials orthogonal with respect to a
measure constructed as a certain combination of the measures $\mu_0$
and $\mu_1$.

\medskip
Let us consider now that the measures $\mu_0$ and $\mu_1$ are
supported on an unbounded interval. There are many asymptotic
results (strong asymptotics) for the monic polynomials $S_{n,
\lambda}$ orthogonal with respect to the inner product
(\ref{producto-lambda}) for a fixed $\lambda$; see for instance
\cite{AMPPR} and \cite{MB} for coherent pairs, \cite{CMM} and
\cite{GLM} for Freud weights and, more recently, the survey
\cite{MM}. But as far as we know, nothing has been said about
asymptotics in the balanced case. In this sense, the first question
that should be answered is: what is the appropriate choice for the
sequence $(\lambda_n)$? We understand by this a sequence of
parameters for which polynomials $S_{n,{\lambda}_n}$ exhibit a
nontrivial asymptotic behavior, depending on both measures $\mu_0$
and $\mu_1$. One of the goals of this paper is to raise that
$\lambda_n = n^{-2}$ is not, in general, the right choice when the
support of $\mu_0$ and $\mu_1$ is unbounded.

The structure of the paper is as follows. In Section 2, we use
heuristic arguments, based on potential theory, about the \lq \lq
size" of ${\lambda}_n$ in order to achieve an appropriate \lq\lq
balancing". In this sense, the Mhaskar--Rakhmanov--Saff numbers turn
out to be a powerful tool. On account of the above results, in
Section 3 we obtain asymptotics for Sobolev polynomials and their
norms for a particular case of Freud weights, which illustrates that
the choice of $\lambda_n$ is accurate.

\section{Selection of the parameters}

We point out some heuristic reasonings about the asymptotic behavior
of the parameters ${\lambda}_n$ in order to balance both terms in
the varying Sobolev inner product $\langle P,Q
\rangle_{{\lambda}_n}$.

Firstly, we recall some basic tools from the classical potential
theory with external field which will be used later on.

Let $\mu$ be a probability measure with support in a closed set
$\Sigma$ of the complex plane. Recall that, the logarithmic
potential $V^{\mu}$ associated with $\mu$ is defined by
$V^{\mu}(z)=-\int \log|z-t|\, d\mu(t)$. Let us assume that
$w(z)=e^{-Q(z)}$ is an admissible and continuous weight function in
$\Sigma$. It is well known that there exists a unique probability
measure $\mu_{w}$, called extremal or equilibrium measure associated
with $w$, minimizing the weighted energy:
$$I_w(\mu)=\int_{\Sigma}\left( V^{\mu}(z)+2Q(z) \right)\,d\mu(z)$$
for all probability measures with support in $\Sigma$. This measure
$\mu_{w}$ is compactly supported and there exists a constant $F_w$
(the modified Robin constant of $\Sigma$) such that $
V^{\mu_{w}}(z)+Q(z)=F_w$ quasi--everywhere on $\supp(\mu_{w})$, see
\cite[Theorem 1.3, p.\ 27]{ST97}. Moreover, if $Q$ is an even
function with some additional properties it can be deduced that
\begin{equation*}
\Vert w^n \, Q_n \Vert_{L_{\infty}(\Sigma)}=\Vert w^n \, Q_n
\Vert_{L_{\infty}(\supp(\mu_{w}))}
\end{equation*}
for every polynomial $Q_n$ of degree $\le n$, see \cite[p.\
203]{ST97}. As a straightforward application of these results, we
can obtain for weighted polynomials a symmetric compact interval on
which its supremum norm lives, more precisely, we have
$$\Vert w \, Q_n \Vert_{L_{\infty}(\Sigma)}=\Vert w \, Q_n
\Vert_{L_{\infty}([-a_n,a_n])}$$ for every polynomial $Q_n$ of
degree $\le n$. The number $a_n$ ($n \ge1$) is the so--called n--th
Mhaskar--Rakhmanov--Saff number for $Q,$ that is, the positive root
of the equation
$$n=\frac{2}{\pi}\, \int_0^1 \frac{a_n\, t\, Q'(a_n\,
t)}{\sqrt{1-t^2}}\, dt\, .$$

The link between the equilibrium measure and the asymptotics of
orthogonal polynomials is given by the following observation: for a
polynomial $Q_n(z)=(z-c_1)(z-c_2) \dots (z-c_n)$ we can write $\log
\vert Q_n(z)\vert=-n\,V^{{\nu}_n}(z)$ where ${\nu}_n$ is the
normalized counting measure on the zeros of $Q_n$, that is,
${\nu}_n=\frac{1}{n}\sum_{i=1}^n \delta_{c_i}$. Then $\vert
w^n(z)\,Q_n(z) \vert ^{1/n}=e^{-(V^{{\nu}_n}(z)+Q(z))}.$

If we denote by $T_{n,w}$ the $n-$th weighted monic Chebyshev
polynomial corresponding to $w$, that is, the solution of the extremal
problem
\begin{equation*}
\inf \{ \Vert w^n\,Q_n \Vert_{L_{\infty}(\Sigma)};\, Q_n(z)=z^n+\dots
\}
\end{equation*}
then
\begin{equation*}
\lim_{n \to \infty} {\Vert w^n\,T_{n,w}
\Vert}_{L_{\infty}(\Sigma)}^{1/n}=e^{-F_w},
\end{equation*}
see \cite[ Theorem 3.1, p.\ 163]{ST97}.

Keeping in mind our balance problem, we are interested in the
asymptotic behavior of the $L_2$--norm in $[-1,1]$ with varying
weights. Since
\begin{equation*}
\lim_{n \to \infty} \left( \frac{\Vert w^n\,Q_n
\Vert_{L_{\infty}([-1,1])}}{\Vert w^n\,Q_n \Vert_{L_{2}([-1,1])}}
\right)^{1/n}=1
\end{equation*}
for every polynomial of degree $n$, (see \cite[Theorem 3.2.1, p.\
65]{ST92}), the asymptotic extremality of ${\Vert w^n\,Q_n
\Vert^{1/n}_{L_{2}([-1,1])}}$ can be thought as the corresponding
one of ${\Vert w^n\,Q_n \Vert^{1/n}_{L_{\infty}([-1,1])}}$. In fact,
if we denote by $ P_{n,w^n}$ the solution of the extremal problem
\begin{equation*}
\inf \{ \Vert w^n\,Q_n \Vert_{L_{2}([-1,1])};\, Q_n(z)=z^n+\dots
\}
\end{equation*}
it can be deduced (see \cite[Theorem 3.3.3, p.\ 78]{ST92}) that
there exists
\begin{equation} \label{pot-lim}
\lim_{n \to \infty} {\Vert w^n\,P_{n,w^n}
\Vert}^{1/n}_{L_{2}([-1,1])}.
\end{equation}

\bigskip

>From now on,   $f_n(x) \sim g_n(x)$ in a domain $D$ will denote that
there are positive constants $C_1$, $C_2$ such that  $C_1 \,g_n(x)
\le f_n(x) \le C_2\,g_n(x)$, for all $x \in D$ and $n$ large enough.

\bigskip

In relation with our problem, we consider the varying Sobolev inner
product $\langle \cdot, \cdot  \rangle_{\lambda_n}$ where
$d\,\mu_i=W^2(x)\,dx$, $i=0,1$. Here, we assume that
$W(x)=e^{-Q(x)}$ is a weight function where $Q: I=(-c,c) \rightarrow
[0, +\infty)$ is a convex, smooth, and even function with
$Q(c^{-})=+\infty=Q((-c)^{+})$ and $Q(x)=0$ only for $x=0$ (we take
$Q$ an even function for simplicity). For these weights $W$, see
\cite[Theorem 4.1, p.\ 95]{LL}, the $L_2$--norm on $I$ for weighted
polynomials is asymptotically equivalent to the $L_2$--norm on a
compact interval. More precisely,
\begin{equation} \label{des-i-mrs}
\Vert W\,Q_n \Vert _{L_2([-a_{n+1}, a_{n+1}])} \le
\Vert W\,Q_n \Vert_{L_2(I)} \le \sqrt{2}\,
 \Vert W\,Q_n \Vert _{L_2([-a_{n+1}, a_{n+1}])}\end{equation}
holds for every $n$ and every polynomial $Q_n$ with degree $\le n$,
where $a_n$ are the Mhaskar--Rakhmanov--Saff numbers associated with
$Q$.

>From (\ref{des-i-mrs}), we deduce that for every polynomial
$Q_n(x)=x^n+ \dots $
 \begin{align} \label{formulalarga}
& \langle Q_n, Q_n \rangle_ {{\lambda}_n } \\
&\sim \int_{-a_{n+1}}^{a_{n+1}} Q_n^2(x)\,W^2(x)\,dx+{\lambda}_n\int_{-a_{n+1}}^{a_{n+1}}
\left( Q_n'(x) \right)^2\,W^2(x)\,dx \notag\\
&=a_{n+1} \left[\int_{-1}^1 Q_n^2(a_{n+1}
t)W^2(a_{n+1} t)dt+{\lambda}_n\int_{-1}^1 \left( Q_n'(a_{n+1} t)
\right)^2W^2(a_{n+1} t)dt \right]\nonumber \\
&=a_{n+1}^{2n+1} \left[ \int_{-1}^1 U_n^2(t)\,W^2(a_{n+1}
t)\,dt+\frac{{\lambda}_n\,n^2}{a_{n+1}^2}\,\int_{-1}^1 V_{n-1}^2(t)
\,W^2(a_{n+1} t)\,dt \right],\notag
\end{align}
where $U_n$ and $V_{n-1}$ are monic polynomials of degree $n$ and
$n-1$, respectively.

Observe that (\ref{formulalarga}) remains true if we take
$d\mu_i=L_i\,W^2(x)\,dx$, $i=0,1$, where $L_0$ and $L_1$ are any
positive constants. At first sight, the presence of the constants
$L_i$ could seem irrelevant but in the next section it will allow us
to give an alternative reading to explain why our selection of
$\lambda_n$ is accurate.

\medskip
Therefore, in order to balance both terms in (\ref{formulalarga}) it
is reasonable to require the following:
\begin{itemize}
\item [i)] $ \lambda_n\, n^2 \sim a_{n+1}^2\, .$
\item [ii)] the asymptotic extremality of the $L_2(W^2(a_{n+1}t),\, [-1,1])$--norm for
monic polynomials of degree $n$ behaves as the corresponding one of
degree $n-1.$
\end{itemize}

The previous results about potential theory lead us to think that a
sufficient condition to get ii) is
\begin{equation} \label{cs-ii}
W^{1/n}(a_{n+1} t)\sim w(t), \quad \forall
t\in(-1,1)
\end{equation}
where $w$ is an admissible and continuous weight function.

\medskip
Concerning the choice of the parameters $\lambda_n$ observe that,
when the support of the measures $\mu_0$ and $\mu_1$ is unbounded,
the size of $\lambda_n$ as $n^{-2}$ is not the right one, in
general. If the weight satisfies (\ref{cs-ii}), the choice of the
parameters depends on the distribution of the measure $W^2(t)\,
dt\,,$ that is, on the corresponding Mhaskar--Rakhmanov--Saff
numbers.

\bigskip

We would like to point out that these ideas can be also applied in a
more general framework. Indeed, consider a Sobolev inner product
with two different weights, $W_0^2$ and $W_1^2$, which are linked in
such a way so that $\langle \cdot , \cdot \rangle_{\lambda_n}$ can
be expressed in terms of only one weight (either $W_0^2$ or $W_1^2$)
satisfying condition (\ref{cs-ii}). Actually, important examples in
this situation are the Hermite coherent pairs. Notice that if the
pair of measures $(W^2_0, W^2_1)$ constitutes a Hermite
symmetrically coherent pair (see \cite{AMPPR} and \cite{MB}), then
either
\begin{align*}
\text{I:} &\quad W_0^2(x)=(x^2+a^2)\, e^{-x^2} \quad \text{and}
\quad W_1^2(x)=e^{-x^2}, \, a\in \R\, ,\quad
\text{or} \\
\text{II:} &\quad W_0^2(x)= e^{-x^2} \quad \text{and} \quad
W_1^2(x)=\frac{e^{-x^2}}{x^2+a^2},\quad a\in \R \setminus\{0\}.
\end{align*}

In both cases we have
\begin{equation*}
\langle Q_n, Q_n \rangle_ {{\lambda}_n } = \int_{\R} \left[
Q_n^2(x)(x^2 + a^2) \right]\,W_1^2(x)\,dx +
{\lambda}_n\int_{\R}\left( Q_n'(x) \right)^2\,W_1^2(x)\,dx,
\end{equation*}
and it is not difficult to check that
\begin{equation*}
\frac{\langle Q_n, Q_n \rangle_ {{\lambda}_n}}{a_{n+2}^{2n+3}} \sim
\int_{-1}^1 U_{n+1}^2(t) W_1^2(a_{n+2}t)\,dt +  \frac{{\lambda}_n
n^2}{a_{n+2}^4}\int_{-1}^1 V_{n-1}^2(t)\,W_1^2(a_{n+2}t)\,dt ,
\end{equation*}
where in each case $a_n$ are the Mhaskar--Rakhmanov--Saff numbers
for the corresponding weight $W_1$, and $U_{n+1}$ and $V_{n-1}$ are
monic polynomials of degree $n+1$ and $n-1$, respectively.

Since $\displaystyle \frac{a_n}{\sqrt{n}} \to \sqrt{2}$, observe
that
\begin{equation*}
\lim_{n \to \infty} W_1^{1/n}(a_{n+2} t) = e^{-t^2}, \quad \forall t\in(-1,1)
\end{equation*}
and therefore, according to the theory stated above, the adequate
choice of $\lambda_n$ should be $\lambda_n \sim a_{n+2}^4 \,
n^{-2}$. In other words, $\lambda_n \sim \text{constant}.$  Hence,
it can be said that the Hermite--Sobolev coherent inner products are
self--balanced.

\section{Freud--Sobolev orthogonal polynomials}

We are going to test the arguments developed in the previous section
for the case of a Sobolev inner product related to Freud weights.
The simplest example corresponds to $W^2_0(x) = W^2_1(x) = e^{-x^2}
$, but this is a trivial case since for any choice of $\lambda_n$
the Sobolev orthogonal polynomial $S_{n,\lambda_n}$ is the $n$--th
monic Hermite polynomial. In this section, we show asymptotics for
the Sobolev orthogonal polynomials with $W^2(x) = W^2_0(x) =
W^2_1(x) = \exp(-x^4).$

Throughout the section, $(P_n)_{n \ge 0}$ denotes the sequence of
monic polynomials orthogonal with respect to the weight $W^2$,
$\Vert \cdot \Vert$ stands  for the $L^2(W^2)$-norm, and
$S_{n,\lambda_n}$ is the monic polynomial which minimizes
\begin{equation*}
\langle Q_n \,, Q_n \rangle_{\lambda_n}=\int_{\R} Q_n^2(x) \, W^2(x) \, dx
+\lambda_n \int_{\R} (Q'_n)^2(x) \, W^2(x) \, dx \,
\end{equation*}
in the class of all monic polynomials of degree $n$.

The Mhaskar--Rakhmanov--Saff numbers for $W(x)
=\displaystyle{\exp(-x^4/2)}$ satisfy $a_n \sim n ^{1/4}$ and
therefore condition (\ref{cs-ii}) holds for $W$. As we have
explained in Section 2, to balance this Sobolev inner product we
must take $\lambda_n \, n^2 \sim a_{n+1}^2,$ that is, $\lambda_n$
like $n^{-3/2}$ when $n \to \infty$.

Next, we study the asymptotic behavior of the ratio $\displaystyle
\frac{S_{n,\lambda_n}}{P_n}$ showing that the choice of $\lambda_n$
provides the reasonable one in a sense we will explain later. For
technical reasons some additional constraints should be imposed on
parameters $\lambda_n ,$ so we deal with a decreasing sequence
$(\lambda_n)$ of positive real numbers such that
\begin{equation} \label{landaajust}
\lim_{n \to \infty} n^{3/2} \lambda_n =L \in [0, +\infty] \,,
\end{equation}
and
\begin{equation} \label{restriccionlandas}
\lim_{n \to \infty} n^{7/4}(\lambda_{n-2} - \lambda_n) = 0 = \lim_{n
\to \infty}
 n^{1/4}\left( \frac{\lambda_{n-2}}{\lambda_n} - 1 \right) \,.
\end{equation}
Notice that the sequence $\lambda_n = n^{-3/2}$ satisfies
(\ref{landaajust}) and (\ref{restriccionlandas}).

\begin{proposition} \label{normas}
Let $(\lambda_n)$ be a decreasing sequence of positive real numbers
which satisfies $\frac{\lambda_{n-2}}{\lambda_n} \to 1$ and $n^{3/2}
\lambda_n \to L \in [0, +\infty].$ Then
\begin{align} \label{limnormas}
\kappa(L):=\lim_{n\to \infty} \frac{\langle S_{n,\lambda_n},
S_{n,\lambda_n} \rangle_{\lambda_n}}{\Vert P_n \Vert^2}  &= \begin{cases}
                1 & \text{if }\, L = 0 \\
                                \frac{2L}{\sqrt{3}} \, \varphi{\left( \frac{20\,L\,+3\,\sqrt{3}}{12\,L}
\,\right)} & \text{if }\, 0< L < +\infty \\
                                                                +\infty &  \text{if }\, L = +\infty\,,
                \end{cases}
\end{align}
where $\varphi (x)=x+\sqrt{x^2-1}$.
\end{proposition}
{\bf Proof.} We consider the Fourier expansion of the polynomial
$P_n$ in terms of the basis $(S_{m,\lambda_n})_{m\ge0}.$ Because of
the weight $e^{-x^4}$ is a symmetric function we have
\begin{equation*}
P_n(z) = S_{n,\lambda_n}(z)
+ \sum_{j=0}^{n-2}\alpha_{j}(\lambda_n)S_{j,\lambda_n}(z)\,,
\end{equation*}
where
\begin{equation*}
\alpha_{j}(\lambda_n)=\frac{\langle P_n,
S_{j,\lambda_n}\rangle_{\lambda_n}}{\langle S_{j,\lambda_n},
S_{j,\lambda_n} \rangle_{\lambda_n}}=\frac{\lambda_n\int_{\R}
P_n'(x)S'_{j,\lambda_n}(x)e^{-x^4}dx}{\langle S_{j,\lambda_n},
S_{j,\lambda_n}\rangle_{\lambda_n}}, \quad 0 \le j \le n-2.
\end{equation*}

Since the orthogonal polynomials $P_n$ satisfy the following
structure relation, (see \cite{N}),
\begin{equation} \label{derivadap}
P'_n(z) = nP_{n-1}(z) + \frac{4 \Vert P_n \Vert^2}{\Vert P_{n-3}
\Vert^2} P_{n-3}(z),
\end{equation}
the coefficients $\alpha_{j}(\lambda_n)$ vanish for $0 \le j<n-2\,
..$ For $j=n-2$ we get
\begin{equation} \label{alfan}
\alpha_{n-2}(\lambda_n) = \frac{4 \,(n-2) \, \lambda_n \, \Vert
P_n \Vert^2}{\langle S_{n-2,\lambda_n}, S_{n-2,\lambda_n}
\rangle_{\lambda_n}} \, ,
\end{equation}
and therefore
\begin{equation} \label{FS}
P_n(z) = S_{n,\lambda_n}(z) + \alpha_{n-2}(\lambda_n)
S_{n-2,\lambda_n}(z), \quad n \ge 3.
\end{equation}

>From now on,  we will write $\kappa_m(\lambda_n) = \langle
S_{m,\lambda_n}, S_{m,\lambda_n}\rangle_{\lambda_n}$, $n,m \ge 0$.

Now, observe that (\ref{FS}) leads to
\begin{align*}
&\kappa_n(\lambda_n)=\langle P_n-\alpha_{n-2}(\lambda_n)
S_{n-2,\lambda_n},P_n-\alpha_{n-2}(\lambda_n)
S_{n-2,\lambda_n}\rangle_{\lambda_n}\\
&=\int_{\R} \left[ \left(P_n-\alpha_{n-2}(\lambda_n)
S_{n-2,\lambda_n}\right)^2 +\lambda_n
\left(P'_n-\alpha_{n-2}(\lambda_n) S_{n-2,\lambda_n}'\right)^2
\right] e^{-x^4}\, dx.
\end{align*}
Then, using (\ref{derivadap}) and the orthogonality of $P_n$ with
respect to the weight function $e^{-x^4}$, we have:
\begin{align*}
\kappa_n(\lambda_n)&=\Vert P_n\Vert^2+n^2\, \lambda_n \, \Vert
P_{n-1}\Vert^2 -8(n-2)\, \lambda_n\, \alpha_{n-2}(\lambda_n)\,
\Vert P_n\Vert^2 \\
&+ 16\, \lambda_n\,\frac{\Vert P_n\Vert^4}{\Vert
P_{n-3}\Vert^2}+\alpha_{n-2}^2(\lambda_n)\,\kappa_{n-2}(\lambda_n)\, .
\end{align*}
Taking into account the value of $ \alpha_{n-2}(\lambda_n)$ given by
(\ref{alfan}), we get
\begin{equation} \label{formulaka}
\kappa_n(\lambda_n)=\Vert P_n \Vert^2 \left(
B_n(\lambda_n)-A_n(\lambda_n)\frac{\Vert P_{n-2}
\Vert^2}{\kappa_{n-2}(\lambda_{n-2})}\right), \, n \ge 3,
\end{equation}
where
\begin{equation*}
A_n(\lambda_n)=16{\lambda_n}^2\,(n-2)^2\,
\frac{\kappa_{n-2}(\lambda_{n-2})}{\kappa_{n-2}(\lambda_n)}
\frac{\Vert P_n \Vert^2}{\Vert P_{n-2} \Vert^2}
\end{equation*}
\begin{equation*}
B_n(\lambda_n)=1+ \lambda_n\, n^2\,\frac{\Vert P_{n-1}
\Vert^2}{\Vert P_n \Vert^2}+16 \lambda_n \, \frac{\Vert P_n
\Vert^2}{\Vert P_{n-3} \Vert^2}.
\end{equation*}

 Next, we study $\lim_n B_n(\lambda_n)$ and $\lim_n
A_n(\lambda_n).$ First, recall that the polynomials $P_n$ satisfy
(see \cite{N})
\begin{equation} \label{normafreud}
\lim_{n\to \infty} \frac{\sqrt{n}\Vert P_{n-1} \Vert^2}{\Vert P_n
\Vert^2} = 2\sqrt{3}.
\end{equation}

On the other hand, $\displaystyle \lim_{n\to \infty}
\frac{\kappa_{n-2}(\lambda_{n-2})}{\kappa_{n-2}(\lambda_n)}=1.$
Indeed, from the assumptions on $\lambda_n$ and using the extremal
property of the norms of monic orthogonal polynomials, we have
\begin{align*}
\kappa_{n-2}(\lambda_n) &\le \kappa_{n-2}(\lambda_{n-2}) \le {\langle
S_{n-2,\lambda_n}, S_{n-2,\lambda_n} \rangle_{\lambda}}_{n-2}  \\
&= \frac{\lambda_{n-2}}{\lambda_n} \left[
\frac{\lambda_n}{\lambda_{n-2}} \Vert S_{n-2,\lambda_n} \Vert^2 +
\lambda_n \Vert S'_{n-2,\lambda_n} \Vert^2 \right] \le
\frac{\lambda_{n-2}}{\lambda_n} \kappa_{n-2}(\lambda_n).
\nonumber
\end{align*}
Since $\frac{\lambda_{n-2}}{\lambda_n} \to 1,$ it follows
\begin{equation} \label{cocientekapas}
\lim_{n\to \infty}
\frac{\kappa_{n-2}(\lambda_{n-2})}{\kappa_{n-2}(\lambda_n)} = 1.
\end{equation}

Firstly, let us suppose that $0 \le L < +\infty.$ Then from
(\ref{cocientekapas}) and (\ref{normafreud}) we deduce that
\begin{equation}
\lim_{n\to \infty} B_n(\lambda_n)=1+\frac{20}{9}\sqrt{3} L, \quad
\text{and} \quad \lim_{n\to \infty} A_n(\lambda_n)= \frac{4}{3} L^2.
\end{equation}
To obtain (\ref{limnormas}) observe that denoting $s_n =
\kappa_n(\lambda_n) /  \Vert P_n \Vert^2$,  (\ref{formulaka})
becomes
\begin{equation} \label{sn}
s_n = B_n(\lambda_n)-A_n(\lambda_n) \frac{1}{s_{n-2}}.
\end{equation}
Writing (\ref{sn}) for even indices and introducing a new sequence
$(q_n)$ by means of $q_{n+1}=s_{2n}q_n$, the above difference
equation becomes
\begin{equation*}
q_{n+1}-B_{2n}(\lambda_{2n})q_n+A_{2n}(\lambda_{2n})q_{n-1}=0,
\end{equation*}
whose characteristic equation
\begin{equation}\label{ecucaract}
q^2-\left(1+\frac{20}{9}\sqrt{3}L\right)q+\frac{4}{3}L^2=0
\end{equation}
has two simple and real roots with distinct moduli. Thus,
Poincar\'e' s Theorem (see, e.g., \cite{poincare}) assures that
$\displaystyle \frac{q_{n+1}}{q_n}=s_{2n}$ converges to a root of
(\ref{ecucaract}). The extremal property of the norms yields
\begin{equation*}
\kappa_n(\lambda_n) \ge \Vert P_n \Vert^2 + \lambda_n \, n^2 \Vert
P_{n-1} \Vert^2,
\end{equation*}
and therefore, using (\ref{normafreud})
\begin{equation*}
l=\lim_{n \to \infty}s_{2n }\ge 1 + \lim_{n \to \infty} \lambda_{2n}
(2n)^2 \frac{\Vert P_{2n-1} \Vert^2}{\Vert P_{2n} \Vert^2}=1+2
\sqrt{3} L.
\end{equation*}
So, it follows easily that $l=\frac{1}{18} \left[ 9+20\sqrt{3}\,L +
\sqrt{768 L^2+360 \sqrt{3}L+81}\right].$ Notice that, if $L \in (0,
+\infty),$ then $l=\frac{2L}{\sqrt{3}} \, \varphi{\left(
\frac{20L+3\sqrt{3}}{12L} \right)}. $

In a similar way, we also prove that $s_{2n+1}$ converges to $l$. As
a conclusion, there exists $\lim_n s_n=l=\kappa(L)$, and so for $L
\in [0, +\infty)$ the Proposition follows.

\medskip
To finish the proof, let us now assume that $L=+\infty.$ From
(\ref{cocientekapas}) and (\ref{normafreud}) we have
\begin{equation*}
\lim_{n \to \infty} \frac{A_n(\lambda_n)}{\left( \lambda_n \,
n^{3/2} \right)^2}=\frac{4}{3}  \quad \text{and} \quad \lim_{n \to
\infty} \frac{B_n(\lambda_n)}{ \lambda_n \, n^{3/2} }=\frac{20}{9}
\sqrt{3}.
\end{equation*}
Upon applying the same technique as in the case $L < +\infty$ and
replacing $s_n$ by $s_n/(\lambda_n \, n^{3/2}),$ we obtain
\begin{equation} \label{sncasoinfty}
\lim_{n \to \infty} \frac{s_n}{\lambda_n \, n^{3/2}}=\lim_{n \to
\infty}\frac{\kappa_n(\lambda_n)}{\lambda_n \, n^{3/2}\Vert P_n
\Vert^2}=2\sqrt{3}.
\end{equation}
Clearly, $\displaystyle{\frac{\kappa_n(\lambda_n)}{\Vert P_n
\Vert^2} \to +\infty}$ when $n$ tends to infinity and we conclude
our statement. $\Box$

\bigskip
The main result of this section is the following:
\begin{theorem} \label{teo-asi-bal-freud}
Let $(\lambda_n)$ be a decreasing sequence of positive real numbers
such that $\lim_n n^{7/4}(\lambda_{n-2}-\lambda_n)=0=\lim_n n^{1/4}
\left( \frac{\lambda_{n-2}}{\lambda_n}-1 \right).$ If
\begin{equation*}
\lim_{n\to \infty} n^{3/2} \lambda_n =L \in [0, +\infty],
\end{equation*}
then
\begin{align*}
\lim_{n\to \infty} \frac{S_{n,\lambda_n}(z)}{P_n(z)}  &= \begin{cases}
1 & \text{if }\, L = 0  \\
\displaystyle \frac{1}{1 - \left[
\varphi\left(\frac{20\,L\,+3\,\sqrt{3}}{12\,L}
\,\right)\right]^{-1}} &\text{if }\, 0 < L < +\infty \\ 3/2 &
\text{if }\, L = +\infty\,,
\end{cases}
\end{align*}
holds uniformly on compact subsets of $\C \setminus \R$.
\end{theorem}

\noindent{\bf Remarks.} 1. The choice $\lambda_n \equiv
\text{constant},$ which corresponds to a non--balanced inner
product, is a particular case of $L=+\infty$, and then Theorem
\ref{teo-asi-bal-freud} recovers the result already obtained in
\cite{CMM}.

2. When $L \in (0, +\infty)$ the above result has also the following
reading. Write
\begin{equation*}
\langle P , Q \rangle_{\lambda_n}=\int_{\R} P(x)\,Q(x) \, W^2(x) \,
dx +\lambda_n \int_{\R} P'(x)\,Q'(x) \,[L\, W^2(x)] \, dx.
\end{equation*}
 If $\lambda_n= n^{-3/2}(1+o(1))$ then $\displaystyle \lim_{n\to
\infty} \frac{S_{n,\lambda_n}}{P_n}$ depends on $L$, that is, on
the ratio of the weights.

However, for any other choice of $\lambda_n$'s the dependence on $L$
disappears, in particular for $\lambda_n=n^{-2}$ (the right choice
in the bounded case) and for $\lambda_n \equiv \text{constant}$ (the
non--balanced case). This shows that our selection of $ \lambda_n$
is accurate since the asymptotic behavior of Sobolev orthogonal
polynomials $S_{n,\lambda_n}$ depends on both measures.

\bigskip
To prove Theorem \ref{teo-asi-bal-freud} we will use the following
result on the strong asymptotics of $P_n$ which appears in
\cite[Section 3]{lopez/rakhmanov:1988}:
\begin{equation} \label{guillermo}
\lim_{n \to \infty} \frac{P_n(z)}{\Vert P_n \Vert} \frac{D_n(z)}{\varphi^{n+1/2}(z/a_n)} =
\frac{1}{\sqrt{2\pi}}
\end{equation}
uniformly on compact subsets of $\C \setminus \R$. Here, $a_n$ are
the Mhaskar--Rakhmanov--Saff numbers associated with the weight
function $W$, $\varphi (z)=z+\sqrt{z^2-1}$ is the conformal mapping
from $\C \setminus [-1, 1]$ onto the exterior of the unit circle,
and
\begin{equation} \label{funcionszego}
D_n(z) = \exp \left( \frac{\sqrt{z^2-a_n^2}}{2\pi} \int_{-a_n}^{a_n}
\frac{-t^4}{(z-t)\sqrt{a_n^2-t^2}}dt \right), \quad z \in \C
\setminus [-a_n, a_n].
\end{equation}
We would like to remark that, for $z \in \C \setminus [-a_n, a_n]$,
\begin{equation*}
D_n(z) = D\left( \frac{1}{\varphi(z/a_n)}, W_n^2 \right),
\end{equation*}
where $W_n^2$ is the weight function on the unit circle $\T$ defined by
\begin{equation*}
W_n^2(e^{i\theta}) = W^2(a_n \cos \theta),\quad \theta \in [-\pi,
\pi],
\end{equation*}
and
\begin{equation*}
D(w, W_n^2) = \exp \left(\frac{1}{2
\pi}\int_0^{2\pi}\frac{e^{i\theta} +w}{e^{i\theta} -w}\log
W_n(\theta) d\theta \right),\quad \vert w \vert <1.
\end{equation*}

It is well known  that $D(. , W_n^2)$ is holomorphic in the open
unit disk $\D$, belongs to the Hardy space $H^2(\D)$, and satisfies:
\begin{itemize}
\item [(1)] $D(w, W_n^2) \not= 0, \quad \textrm{for} \quad w \in \D$
\item [(2)] $D(0, W_n^2) > 0$
\item [(3)] for almost every $\zeta$ in the unit circle, $D(., W_n^2)$ has nontangential
boundary values $D(\zeta, W_n^2)$ such that $\vert D(\zeta, W_n^2)
\vert^2 = W_n^2(\zeta),$
\end{itemize}
\noindent (see, for instance, \cite{R}).

Next, we prove a technical result that will be also used in the
proof of Theorem \ref{teo-asi-bal-freud}.
\begin{lemma}
Assume that the sequence $(\lambda_n)$ satisfies the same conditions as in Theorem
\ref{teo-asi-bal-freud}, then
\begin{equation*}
\lim_{n \to \infty}
\frac{S_{n,\lambda_{n-2}}(z)-S_{n,\lambda_n}(z)}{P_n(z)}=0,
\end{equation*}
uniformly on compact subsets of $\C \setminus \R$.
\end{lemma}
\noindent {\bf Proof.} On account of (\ref{guillermo}), it suffices
to prove that
\begin{equation*}
\lim_{n \to \infty}
\frac{S_{n,\lambda_{n-2}}(z)-S_{n,\lambda_n}(z)}{\Vert P_n \Vert \,
\varphi^{n+1/2}(z/a_n)}D_n(z)=0
\end{equation*}
uniformly on compact subsets of $\C \setminus \R$. To see this we
will prove:
\begin{itemize}
\item [i)] for every compact set $K$ in $\C \setminus \R$, there exists a
constant $M_K$, not depending on $n$, such that for $n$ large enough
\begin{align*} \nonumber
&\sup_{z \in K} \left| \frac{S_{n,\lambda_{n-2}}(z)-
S_{n,\lambda_n}(z)}{\Vert P_n \Vert \, \varphi^{n+1/2}(z/a_n)}D_n(z) \right|^2 \\
&\le M_K \,a_n \int_{-a_n}^{a_n} \frac{|S_{n,\lambda_{n-2}}(x)-
S_{n,\lambda_n}(x)|^2}{\Vert P_n \Vert^2} \, W^2(x) \, dx,
\end{align*}
and
\item [ii)]
\begin{equation*}
\lim_{n \to \infty} a_n \int_{-a_n}^{a_n}
\frac{|S_{n,\lambda_{n-2}}(x)-S_{n,\lambda_{n}}(x)|^2}{\Vert P_n
\Vert^2} \, W^2(x)\,dx=0.
\end{equation*}
\end{itemize}

The key idea to prove $i)$ is to use the conformal mapping
$\varphi(z/a_n)$ which applies  $\C \setminus [-a_n, a_n]$ onto
$\Omega = \{ z \in \C; \vert z \vert > 1 \}$, and the Cauchy
integral representation for functions in $H^2(\Omega)$. Here,
$H^2(\Omega)$ denotes the space of analytic functions $f$ in
$\Omega$, with limit at $\infty$ and such that $f(\frac{1}{z})$
belongs to the Hardy space $H^2(\D)$. From the Cauchy integral
representation for functions in $H^2(\D)$, see \cite{R}, we have
that if $f \in H^2(\Omega)$ then
\begin{equation} \label{integralcauchy}
f(w)= - \frac{1}{2 \pi i}\int_{\vert \zeta \vert =1}
\frac{f^{*}(\zeta)}{\zeta - w} \, \frac{w}{\zeta} \, d \zeta,\quad w
\in \Omega
\end{equation}
where $f^{*}(\zeta) = \lim_{r \searrow 1} f(r \, \zeta)$ and the
unit circle is positively oriented.

In order to prove $i)$, given a compact set $K$ in $\C \setminus
\R$, there exists an absolute constant $C_K > 0$ such that
\begin{equation*}
 \vert \sqrt{z^2-a_n^2} \vert \ge C_K, \quad \forall z \in K, \quad \forall n
 \ge 0.
\end{equation*}
Therefore, if $z \in K$,
\begin{align*}
&\left| \frac{S_{n,\lambda_{n-2}}(z) - S_{n,\lambda_n}(z)}{\Vert P_n \Vert
\varphi^{n+1/2}(z/a_n)}  D\left(\frac{1}{\varphi(z/a_n)}, W^2_n \right) \right|^2 \\
&\le \frac{1}{C_K} \left| \frac{S_{n,\lambda_{n-2}}(z) - S_{n,\lambda_n}(z)}{\Vert P_n \Vert
 \varphi^{n+1/2}(z/a_n)}  D\left(\frac{1}{\varphi(z/a_n)}, W^2_n \right) \right|^2 \, \vert
\sqrt{z^2 - a_n^2} \vert \\ &= \frac{1}{C_K} \vert F_n(w) \vert
\end{align*}
where
\begin{equation*}
F_n(w) = \left[ \frac{(S_{n,\lambda_{n-2}}-S_{n,\lambda_n})(a_n \,
\varphi^{-1}(w))}{\Vert P_n \Vert  w^{n+1/2}}D\left(\frac{1}{w},
W_n^2 \right) \right]^2 a_n \sqrt{(\varphi^{-1}(w))^2-1},
\end{equation*}
with $w =\varphi(z/a_n)$.

It is easy to check that $F_n \in H^2(\Omega)$ and its boundary
values are
\begin{equation*}
F_n^{*}(e^{i \theta})=\frac{(S_{n,\lambda_{n-2}}-
S_{n,\lambda_n})^2(a_n  \cos \theta)}{\Vert P_n \Vert^2  e^{i
(2n+1)\theta}}W^2(a_n \, \cos \theta)a_n \sqrt{\cos^2 \theta - 1}.
\end{equation*}

Moreover, if we denote by $K_n = \{ \varphi(z/a_n) ; z \in K \}$,
straightforward computations yield that there exists an absolute
constant $A_K > 0$ such that the distance between $K_n$ and the unit
circle satisfies $d(K_n, \T) \ge A_K/a_n$ for $n$ large enough.
Then, from the integral formula (\ref{integralcauchy}) applied to
$F_n$ we have for $w \in K_n$
\begin{align*}
\vert F_n(w) \vert &\le B_K \,a_n \int_{\vert \zeta \vert=1} \vert
F_n^{*} (\zeta) \vert \,
\vert d\, \zeta \vert \\
&= B_K \, a_n \int_{-\pi}^{\pi} \frac{(S_{n,\lambda_{n-2}}-
S_{n,\lambda_{n}})^2(a_n \, \cos\theta)}{\Vert P_n \Vert^2} \,
W^2(a_n \, \cos\theta) \, a_n
\vert \sen \, \theta \vert \, d \theta \\
&= 2 \, B_K \, a_n \,  \int_{-a_n}^{a_n} \frac{(S_{n,\lambda_{n-2}}
- S_{n,\lambda_{n}})^2(x)}{\Vert P_n \Vert^2} \,W^2(x) \,d x,
\end{align*}
where $B_K$ is an absolute positive constant depending only on $K$. So $i)$
is proved.

\medskip
In order to deduce $ii)$, observe that
\begin{align*}
&\int_{\R} |S_{n,\lambda_{n-2}}(x)-S_{n,\lambda_n}(x) |^2 W^2(x)\,
dx \le \langle S_{n,\lambda_{n-2}}-S_{n,\lambda_n},
S_{n,\lambda_{n-2}}-S_{n,\lambda_n} \rangle_{\lambda_n} \\
&=\langle S_{n,\lambda_{n-2}}, S_{n,\lambda_{n-2}}
\rangle_{\lambda_n} -
\langle S_{n,\lambda_n}, S_{n,\lambda_n} \rangle_{\lambda_n} \\
&= \kappa_n(\lambda_{n-2})+(\lambda_n - \lambda_{n-2})\int_{\R}
|S'_{n,\lambda_{n-2}}(x)|^2 \, W^2(x) \,dx - \kappa_n(\lambda_n) \\
&\le \kappa_n(\lambda_{n-2}) - \kappa_n(\lambda_n).
\end{align*}
Therefore, for every $n$  we get
\begin{equation*}
a_n \int_{-a_n}^{a_n} \frac{|S_{n,\lambda_{n-2}}(x) -
S_{n,\lambda_{n}}(x)|^2}{\Vert P_n \Vert^2} W^2(x) dx \le
 a_n \, \frac{\kappa_n(\lambda_{n-2}) - \kappa_n(\lambda_n)}{\Vert P_n \Vert^2}.
\end{equation*}
Finally, since $a_n \sim n^{1/4}$ it is enough to prove that
$\displaystyle n^{1/4} \, \frac{\kappa_n(\lambda_{n-2}) -
\kappa_n(\lambda_n)}{\Vert P_n \Vert^2}$ tends to $0$ when $n$ tends
to infinity. Indeed, since
\begin{equation*}
0 \le \kappa_n(\lambda_{n-2})- \kappa_n(\lambda_n) \le \left(1-
\frac{\lambda_n}{\lambda_{n-2}}\right) \kappa_n(\lambda_{n-2}),
\end{equation*}
we have
\begin{equation*}
\displaystyle n^{1/4} \, \frac{\kappa_n(\lambda_{n-2})-
\kappa_n(\lambda_n)}{\Vert P_n \Vert^2} \le n^{1/4} \, \left(1-
\frac{\lambda_n}{\lambda_{n-2}}\right)
\frac{\kappa_n(\lambda_n)}{\Vert P_n \Vert^2}
\frac{\kappa_n(\lambda_{n-2})}{\kappa_n(\lambda_n)}.
\end{equation*}
Now, taking into account that $\lim_n
\frac{\kappa_n(\lambda_{n-2})}{\kappa_n(\lambda_n)}= 1$, it suffices
to keep in mind Proposition \ref{normas}, (\ref{restriccionlandas}),
and (\ref{sncasoinfty}) to conclude $ii)$ and therefore the proof of
the Lemma. $\Box$

\noindent {\bf Proof of Theorem \ref{teo-asi-bal-freud}.} The
algebraic relation between the polynomials $P_n$ and the Sobolev
polynomials given by (\ref{FS}) can be rewritten for $\lambda_{n-2}$
as
\begin{align*}
P_n(z)&=S_{n,\lambda_{n-2}}(z)+\alpha_{n-2}(\lambda_{n-2})\,
S_{n-2,\lambda_{n-2}}(z)\\
&=S_{n,\lambda_{n}}(z)+S_{n,\lambda_{n-2}}(z)-S_{n,\lambda_{n}}(z)+\alpha_{n-2}(\lambda_{n-2})\,
S_{n-2,\lambda_{n-2}}(z).
\end{align*}
Then, dividing both hand sides of the above expression by $P_n(z)$,
we obtain
\begin{equation} \label{efes}
f_n(z) = b_n(z)f_{n-2}(z)+c_n(z), \quad z \in {\C} \setminus \R
\end{equation}
where
\begin{align*}
f_n(z) &= \frac{S_{n,\lambda_n}(z)}{P_n(z)}, \quad
b_n(z)=-\alpha_{n-2}(\lambda_{n-2})\,\frac{P_{n-2}(z)}{P_n(z)},\\
c_n(z) &= 1-
\frac{S_{n,\lambda_{n-2}}(z)-S_{n,\lambda_n}(z)}{P_n(z)}.
\end{align*}

Firstly, we study the limits of the sequences $(b_n(z))$ and
$(c_n(z))$. As a consequence of Lemma 1 we know that
\begin{equation*}
\lim_{n \to \infty} c_n(z)=1
\end{equation*}
uniformly on compact subsets of $\C \setminus \R$.

With regard to $(b_n(z))$, if $L \in [0,+\infty)$ from Proposition
\ref{normas} and (\ref{normafreud})
\begin{equation*}
\frac{\alpha_{n-2}(\lambda_{n-2})}{\sqrt{n-2}}=4 \lambda_{n-2}
(n-2)^{3/2} \frac{\Vert P_{n-2}
\Vert^2}{\kappa_{n-2}(\lambda_{n-2})} \frac{\Vert P_n
\Vert^2}{(n-2)\,\Vert P_{n-2} \Vert^2} \to \frac{L}{3\kappa(L)}.
\end{equation*}
Moreover, for the monic polynomials $P_n$ it is
known, (see \cite{lopez/rakhmanov:1988}), that
\begin{equation*}
\lim_{n\to \infty}\frac{\sqrt{n-2}\,P_{n-2}(z)}{P_n(z)}=
-2\,\sqrt{3},
\end{equation*}
uniformly  on compact subsets of $\C \setminus \R$. Both results
lead to
\begin{equation*}
\lim_{n\to \infty} b_n(z)= \frac{2L}{\sqrt{3}\,\kappa(L)}=
\begin{cases} 0 & \, \text{if }
  L = 0 \\
\frac{1}{\varphi \left( \frac{20\,L\,+3\,\sqrt{3}}{12\,L} \right)} &
\,  \text{if } 0< L < +\infty \\
\end{cases}
\end{equation*}
uniformly on compact subsets of $\C \setminus \R$. In the case $L =
+\infty$, using formula (\ref{sncasoinfty}) we get $\lim_n b_n(z)=
1/3$, uniformly on compact subsets of $\C \setminus \R$.

Finally, observe that the functions $f_n$, $b_n$ and $c_n$  are
analytic in $\C \setminus \R$. Since for $L \in [0, +\infty]$ we
have $\lim_n b_n(z) = b_L$, with $|b_L| < 1,$ and $\lim_n c_n(z) =
1$ uniformly on compact subsets of $\C \setminus \R$, we can deduce
that
\begin{equation*}
\lim_{n\to \infty} f_n(z)= \frac{1}{1-b_L}
\end{equation*}
uniformly on compact subsets of $\C \setminus \R$. Indeed, for a
fixed compact set $K \subset  \C \setminus \R$ there exist constants
$r \in (0,1), \, R > 1$ and a positive integer number $n_0$ such
that
\begin{equation*}
\vert b_n(z) \vert \le r, \quad \vert c_n(z) \vert \le R , \quad
\text{for} \quad n \ge n_0, \quad z \in K.
\end{equation*}
Thus
\begin{equation*}
\vert f_n(z) \vert \le r \vert f_{n-2}(z) \vert + R , \quad
\text{for}  \quad n \ge n_0 ,\quad z \in K,
\end{equation*}
and therefore we deduce that the sequence $(f_n)$ is uniformly
bounded on compact subsets of $\C \setminus \R$.

>From (\ref{efes}), we can write
\begin{equation*}
f_n(z)- \frac{1}{1-b_L}=b_L \left[ f_{n-2}(z)- \frac{1}{1-b_L}
\right] + \varepsilon_n(z),
\end{equation*}
with
\begin{equation*}
\varepsilon_n(z)=(b_n(z) - b_L) f_{n-2}(z)+c_n(z)-1.
\end{equation*}
Notice that $\lim_n \varepsilon_n(z)=0$, uniformly on compact
subsets of $\C \setminus \R$. From the fact $|b_L| < 1$, it is easy
to deduce that
\begin{equation*}
\lim_{n\to \infty} f_n(z)=\frac{1}{1-b_L},
\end{equation*}
uniformly on compact subsets of $\C \setminus \R$. Taking into
account the value of $b_L$ with $L \in [0, +\infty]$, the Theorem is
proved.  $\Box$

\subsection*{Acknowledgements}
We express our gratitude to Professor Andrei Mart\'{\i}nez--Finkelshtein
for his suggestions and fine comments on this research. We also
thank the referee for his observations and comments which have
improved this article.

\end{document}